\documentclass[a4paper,11pt]{article}
\usepackage{amsmath, amsfonts, amssymb}
\usepackage[english, russian]{babel}
\begin{document}
\makeatletter
\makeatother
\small \begin{center}{\textit{ In the name of
Allah, the Beneficent, the Merciful.}}\end{center}
\large
\begin{center}
\textbf{On classification of finite dimensional algebras}\\
\textbf{Ural Bekbaev}\\
{Department of Science in Engineering, Faculty of Engineering,\\
  IIUM, Kuala Lumpur, Malaysia;}\\
{bekbaev@iium.edu.my}
\end{center}
\small \begin{center} Abstract\end{center}
Classification and invariants, with respect to basis changes, of finite dimensional algebras are considered. An invariant open, dense (in the Zariski topology) subset of the space of structural constants is defined. The algebras with structural constants from this set are classified and a basis to the field of invariant rational functions of structural constants is provided.
\large
\section{Introduction}

The classification of finite dimensional algebras is an important problem in Algebra. For example, the classification of finite dimensional simple and semi-simple associative algebras by
Wedderburn, the classification of finite dimensional simple and semi-simple Lie algebras by Cartan
are considered as key results in the theory of corresponding algebras. In these both and many other cases the used approach to the classification problem is structural (basis free, invariant). Unlike the structural approach in this paper we are going to deal with all algebras by the use of their structural constants. For similar approach in small dimensional case one can consider [1-3]. In general case similar approach is considered in [4]. We also consider the
equivalence problem of algebras in general case. Though there are some intersecting results in this paper with [4] our approach is more simple and more
constructive than of [4].  For the given dimension we show how to construct an invariant, open, dense subset of the space of structural constants and classify all algebras who's system of structural constants are in this set. We provide a  basis for the field of invariant rational functions of structural constants as well.

The paper is organized in the following way. The key results which are used to obtain the classification and invariants of algebras are presented in Section 2. Section 3 is a realization of Section 2 results in the case of representation of general linear group in the space structural constants.

\section{Preliminaries}

In this section we consider a linear representation of a subgroup of the general linear group and under an assumption prove some general results about the equivalence and invariance problems with respect to this subgroup.

  Let $n$, $m$  be any natural numbers, $\tau: (G,V)\rightarrow V$ be a fixed linear algebraic representation of an algebraic subgroup $G$ of $GL(m,F)$ over $V$, where $F$ is any  field and $V$ is $n$-dimensional
 vector space over $F$. Further we consider this representation under the following assumption:

 {\bf Assumption.}There exists a nonempty $G$-invariant subset $V_0$ of $V$ and an algebraic map $P: V_0\rightarrow G$ such that
 $P(\tau( g,\mathbf{v}))=P(\mathbf{v}) g^{-1}$ whenever  $\mathbf{v}\in V_0$ and $g\in G$.

{\bf Theorem 2.1.} Elements $\mathbf{u},\mathbf{v}\in V_0$ are $G$-equivalent, that is $\mathbf{u}=\tau( g,\mathbf{v})$ for some $g\in G$, if and only if $\tau( P(\mathbf{u}),\mathbf{u})=\tau( P(\mathbf{v}),\mathbf{v})$.

{\bf Proof.} If $\mathbf{u}=\tau( g,\mathbf{v})$ then $\tau( P(\mathbf{u}),\mathbf{u})=\tau( P(\tau( g,\mathbf{v})),\tau( g,\mathbf{v}))=$ \[\tau( P(\mathbf{v}) g^{-1},\tau( g,\mathbf{v}))=\tau( P(\mathbf{v}),\tau( g^{-1},\tau( g,\mathbf{v})))=\tau( P(\mathbf{v}),\mathbf{v}).\]
Visa versa, if $\tau( P(\mathbf{u}),\mathbf{u})=\tau( P(\mathbf{v}),\mathbf{v})$ then \[ \tau( P(\mathbf{u})^{-1}P(\mathbf{v}),\mathbf{v})=\tau(( P(\mathbf{u}))^{-1},\tau( P(\mathbf{v}),\mathbf{v}))= \tau(( P(\mathbf{u}))^{-1},\tau( P(\mathbf{u}),\mathbf{u}))=\mathbf{u}\] that is $\mathbf{u}=\tau( g,\mathbf{v})$, where  $g=P(\mathbf{u})^{-1}P(\mathbf{v})$.

This proposition shows that the system of components of $\tau( P(\mathbf{x}),\mathbf{x})$ is a separating system of invariants for the $G$-orbits in $V_0$, where $\mathbf{x}=(x_1, x_2,...,x_n)$ is an algebraic independent system of variables over $F$.

Further in this paper it is assumed that $F$ is algebraically closed field, $V_0$ in the above Assumption is an open, dense (in Zariski topology) $G$-invariant subset of $V$. In such case for any $\mathbf{u},\mathbf{v}\in V_0$ one has $P(\tau( P(\mathbf{u}),\mathbf{v}))=P(\mathbf{v}) P(\mathbf{u})^{-1}$ and due to density of $V_0$ in $V$ one has
\[P(\tau( P(\mathbf{y}),\mathbf{x}))=P(\mathbf{x}) P(\mathbf{y})^{-1},\] where $\mathbf{y}=(y_1, y_2,...,y_n)$ is also an algebraic independent system of variables over $F$.

{\bf Theorem 2.2.} The field of $G$-invariant rational functions $F(\mathbf{x})^G$ is generated over $F$ by the system of components of $\tau( P(\mathbf{x}),\mathbf{x})$.

{\bf Proof.}  It is evident that all components of  $\tau( P(\mathbf{x}),\mathbf{x})$ are in $ F(\mathbf{x})^G$. If $f(\mathbf{x})=f(\tau( g,\mathbf{x}))$ for all $g\in G$ then, in particular, $f(\mathbf{x})=f(\tau( P( \mathbf{u}),\mathbf{x}))$ whenever $\mathbf{\mathbf{\mathbf{u}}}\in V_0$. It implies, as far as $V_0$ is dense in $V$, that for the variable vector $\mathbf{y}$ the equality \[f(\mathbf{x})=f(\tau( P( \mathbf{y}),\mathbf{x}))\]
 holds true. In $\mathbf{y}=\mathbf{x}$ case one gets that $f(\mathbf{x})=f(\tau( P( \mathbf{x}),\mathbf{x}))$.

 {\bf Remark 2.1.} In the previous variant of this paper the pure transcendence of the extension $F\subset F(\mathbf{x})^{GL(m,F)}$ is stated with some justification. It turns out  that the
  justification is not valid and therefore it is still not clear if the assumption alone is sufficient to state the pure transcendence of the extension $F\subset F(\mathbf{x})^G$. Let $V_{00}$ stand for $\{v\in V_0: P( v)=I_m\}$, where $I_m$ stands for the $m$-order identity matrix. If $w=\tau( g,v)$ and $w,v\in V_{00}$ then  $I_m=P( w)=P(\tau( g,v))=P( v) g^{-1}$ that is $g=I_m$. It implies that each $G$-orbit from $V_0$ has only one common element with $V_{00}$. Each $G$-invariant function on $V_0$ is uniquely defined by its values on $V_{00}$.
Note that the sets $G\times V_{00}$ and $V_0$ are bi-rationally isomorph. Indeed the map
  \[Q: G\times V_{00}\rightarrow V_0,\ \mbox{where}\ Q( g, \mathbf{v})= \tau( g,\mathbf{v}) \] is invertible with
$Q^{-1}: V_0\rightarrow G\times V_{00},\ \mbox{where}\ Q^{-1}( \mathbf{v})= (P(\mathbf{v})^{-1},\tau( P(\mathbf{v}),\mathbf{v}))$.

 {\bf Corollary 2.1.} The field $F(\mathbf{x})$ is generated over $F(\mathbf{x})^G$ by the system of components of $P(\mathbf{x})$.

{\bf Proof.} Indeed $F(\mathbf{x})^G( P(\mathbf{x}))=F(\tau( P(\mathbf{x}),\mathbf{x}))( P(\mathbf{x}))=F(\tau( P(\mathbf{x}),\mathbf{x}),P(\mathbf{x}))$ and \\ $\tau( P(\mathbf{x})^{-1},\tau( P(\mathbf{x}),\mathbf{x})) =\mathbf{x}$ and therefore $F(\mathbf{x})^G(  P(\mathbf{x}))=F(\mathbf{x})$.

{\bf Proposition 2.1.} The equality $\mbox{trdeg} F( P(\mathbf{x}))/F= \mbox{dim} G$ holds true.

{\bf Proof.} To prove the equality it is enough to show equality of the vanishing ideals of $P(\mathbf{x})$ and $G$.
If polynomial $p$ vanished on $P(\mathbf{x})$, that is $p[P(\mathbf{x})]=0$, then  $p[P(\tau( g,\mathbf{x},))]=p[P(\mathbf{x}) g^{-1}]=0$. In particular, $p[g]=0$ for any $g\in G$ that is $p$ vanishes on $G$ as well.

If $p[g]=0$ for any $g\in G$ then, in particular, $p[P(\mathbf{u})]=0$ for any $\mathbf{u}\in V_0$. Due to density of $V_0$ in $V$ one has
$p[P(\mathbf{x})]=0$.

{\bf Theorem 2.3.} The equality $\mbox{trdeg}F( \mathbf{x})^G/F=n-\mbox{dim}G$ holds true.

{\bf Proof.} A proof of the theorem can be deduced from Remark 1. Here is another more detailed proof. Let $\widetilde{P(\mathbf{x})}$ stand for any system of entries of $P(\mathbf{x})$ which is a transcendence basis for the field $F( P(\mathbf{x}))$ over $F$. We show that the system $\widetilde{P(\mathbf{x})}$ is algebraic independent over $F(\mathbf{x})^G$ as well. Indeed let $p[\widetilde{(y_{ij})_{i,j=1,2,...,m}}]$ be any polynomial over $F(\mathbf{x})^G$ for which  $p[ \widetilde{P(\mathbf{x})}]=0$ that is $p_\mathbf{v}[\widetilde{P(\mathbf{v})}]=0$ for all $\mathbf{v}\in V_1$, where $V_1$ is a $G$-invariant nonempty open subset of $V_0$, where $p_\mathbf{v}[\widetilde{(y_{ij})_{i,j=1,2,...,m}}]$ stands for the polynomial obtained from $p[\widetilde{(y_{ij})_{i,j=1,2,...,m}}]$ by substitution $\mathbf{v}$ for $\mathbf{x}$.
The equality  $0=p_\mathbf{v}[\widetilde{P(\mathbf{v})}]=p_{\tau( g,\mathbf{v})}[\widetilde{P(\tau( g,\mathbf{v}))}]=p_\mathbf{v}[\widetilde{P(\mathbf{v}) g^{-1}}]$ implies that $p_\mathbf{v}[\widetilde{g}]=0$ for any $g\in G$. Therefore $p_\mathbf{v}[\widetilde{P(\mathbf{x})}]=0$, that is $p_\mathbf{v}[\widetilde{(y_{ij})_{i,j=1,2,...,m}}]$ is zero polynomial for any $\mathbf{v}\in V_1$. It means that $p[\widetilde{(y_{ij})_{i,j=1,2,...,m}}]$ is zero polynomial itself. Now due to $F\subset F(\mathbf{x})^G\subset F(\mathbf{x})$, $\mbox{tr.deg.}F(\mathbf{x})/F=n$ and Corollary 2.1 one has the required result.

{\bf Corollary 2.2} The transcendence degree of $F(\mathbf{x})^{GL(m,F)}$ over $F$ equals to $n-m^2$ and the field extension $F(\mathbf{x})^{GL(m,F)}\subset F(\mathbf{x})$ is a pure transcendental extension. For a transcendental basis one can take the system of components of $P(\mathbf{x})$.

{\bf Question 2.1.} Under the assumption for $G=GL(m,F)$ is it true that $F\subset F(\mathbf{x})^{GL(m,F)}$ is also a pure transcendental extension?

{\bf Question 2.2.} Is it a typical situation for any representation of $G=GL(m,F)$ that if one of the extensions $F\subset F(\mathbf{x})^G$, $F(\mathbf{x})^G\subset F(\mathbf{x})$ is a pure transcendental extension then the second one also  has the same property?

\section{Classification of algebras }
\subsection{General case }

In this paper we use the standard notation (the Einstein notation) for tensors as well as the matrix representation for tensors which is more convenient in dealing with equivalence and invariance problems of tensors with respect to basis changes. The use of matrix representation for tensors makes the descriptions more transparent as well.

 Let us consider any $m$ dimensional algebra $W$ with multiplication $\cdot$ given by a bilinear map $(\mathbf{u},\mathbf{v})\mapsto \mathbf{u}\cdot \mathbf{v}$. If $e=(e^1,e^2,...,e^m)$ is a basis for $W$ then one can represent the bilinear map by a matrix $A\in Mat(m\times m^2;F)$ such that \[\mathbf{u}\cdot \mathbf{v}=eA(u\otimes v)\] for any $\mathbf{u}=eu,\mathbf{v}=ev,$ where $u=(u_1, u_2,..., u_m), v=(v_1, v_2,..., v_m)$ are column vectors.
So the binary operation (bilinear map, tensor) is presented by the matrix $A\in Mat(m\times m^2;F)$ with respect to the basis $e$.
Further we deal only with such matrices of rank $m$.

 If $e'=(e'^1,e'^2,...,e'^m)$ is also a basis for $W$, $g\in G=GL(m,F)$, $e'g=e$ and  $\mathbf{u}\cdot \mathbf{v}=e'B(u'\otimes v')$, where $\mathbf{u}=e'u',\mathbf{v}=e'v'$, then $\mathbf{u}\cdot \mathbf{v}=eA(u\otimes v)=e'B(u'\otimes v')=eg^{-1}B(gu\otimes gv)=eg^{-1}B(g\otimes g)(u\otimes v)$ as far as $\mathbf{u}=eu=e'u'=eg^{-1}u',\mathbf{v}=ev=e'v'=eg^{-1}v'$. Therefore the equality
\begin{equation}B=gA(g^{-1})^{\otimes 2}\end{equation} is valid.

Now let $\tau$ stand for the representation of $G=GL(m,F)$ on the $n=m^3$ dimensional vector space $V=Mat(m\times m^2;F)$ defined by \[ \tau: (g,A)\mapsto B=gA(g^{-1}\otimes g^{-1}).\]

To have Theorems 2.1-2.3 for this case we will construct a map $P: V_0\rightarrow GL(m,F)$ with property (1) in the following way. For any natural number $k$ due to $(2)$ one has \begin{equation}B^{\otimes k}=g^{\otimes k}A^{\otimes k}(g^{-1})^{\otimes 2k}\end{equation}
Let us consider all its possible contractions with respect to $k$ upper and $k$ lower indices. It is clear that the result of each of such contraction will be  $f(B)=f(A)(g^{-1})^{\otimes k}$ type equality, where $f(A)$ is a row vector with $m^k$ entries.

In $k=1$ case one gets the following $2^11!=2$ different row equalities: $\mathbf{Tr_1}(B)=\mathbf{Tr_1}(A)g^{-1}, \ \mathbf{Tr_2}(B)=\mathbf{Tr_2}(A)g^{-1},$ where $\mathbf{Tr_1}(A)$ stands for the row vector  with entries $A^j_{j,i}=\sum_{j=1}^nA^j_{j,i}$- the contraction on the first upper and lower indices and $\mathbf{Tr_2}(A)$ stands for the row vector  with entries $A^j_{i,j}=\sum_{j=1}^nA^j_{i,j}$- the contraction on the first upper and second lower indices.

In $k=2$ case one gets the following $2^22!+2^11!=10$ different row equalities:
\[\mathbf{Tr_i}(B)\otimes \mathbf{Tr_j}(B)=(\mathbf{Tr_i}(A)\otimes \mathbf{Tr_j}(A))(g^{-1})^{\otimes 2}, \ \mathbf{Tr_i}(B)B=\mathbf{Tr_i}(A)A(g^{-1})^{\otimes 2},\] where $i,j=1,2,$
and \[(B^i_{j,p}B^j_{i,q})=(A^i_{j,p}A^j_{i,q})(g^{-1})^{\otimes 2},\ (B^i_{j,p}B^j_{q,i})=(A^i_{j,p}A^j_{q,i})(g^{-1})^{\otimes 2},\] \[(B^i_{p,j}B^j_{i,q})=(A^i_{p,j}A^j_{i,q})(g^{-1})^{\otimes 2},\ (B^i_{p,j}B^j_{q,i})=(A^i_{p,j}A^j_{q,i})(g^{-1})^{\otimes 2}.\]

In any $k$ case only the number of contractions of $A^{\otimes k}$ when all $k$ different  upper indices are contracted with lower indices of different $A$ is \[(2k)\times (2(k-1))\times (2(k-2))\times...\times 2=
   2^kk!.\] In general it is nearly clear that the corresponding resulting system of $2^kk!$ rows depending on the variable matrix $A:=X=(X^i_{j,k})_{i,j,k=1,2,...,m}$ is linear independent over $F$. But for big enough $k$ the inequality $2^kk!\geq m^k$ holds true as well. Therefore in general for big enough $k$ it is possible to choose $m^k$ contractions (rows) among the all contractions of $X^{\otimes k}$ for which the matrix $Q(X)$
   consisting of these $m^k$ rows is a nonsingular matrix. For the matrix $Q(X)$ one has equality $Q(Y)=Q(X)(g^{-1})^{\otimes k}$ whenever $g\in G$, $Y=gX(g^{-1})^{\otimes 2}$.

  Now note that for any $A\in \{X: \det(Q(X))\neq 0\}$ and $g\in G$ one has, for example,
$(B\otimes (\mathbf{Tr_1}(B))^{\otimes k-2})Q(B)^{-1}=$ \[g(A\otimes (\mathbf{Tr_1}(A))^{\otimes k-2})(g^{-1})^{\otimes k}(Q(A)(g^{-1})^{\otimes k})^{-1}=g(A\otimes (\mathbf{Tr_1}(A))^{\otimes k-2})Q(A)^{-1}.\] Therefore if $P(A)^{-1}$ stands for arbitrary  nonsingular $m\times m$ size sub-matrix of $(A\otimes (\mathbf{Tr_1}(A))^{\otimes k-2})Q(A)^{-1}$ then one has the  equality $P(B)^{-1}=gP(A)^{-1}$, where $g\in G$, $B=gA(g^{-1})^{\otimes 2}$. It implies that whenever $A\in V_0=\{A: \det(P(A))\det(Q(A))\neq 0\}$ the equality
$P(B)=P(A)g^{-1}$ holds true for any $g\in G$ and $B=gA(g^{-1})^{\otimes 2}$. Note that \[V_0=\{A: \det(P(A))\det(Q(A))\neq 0\}\] is a $G$-invariant, open and dense subset of $V$.

Therefore we have the following results.

{\bf Theorem 3.1.} Two algebras with matrices of structural constants $A,B\in V_0$ are same (isomorph) algebras if and only if
 \[P(A)A(P(A)^{-1}\otimes P(A)^{-1})=P(B)B(P(B)^{-1}\otimes P(B)^{-1}).\]

{\bf Theorem 3.2.} The field of $G$-invariant rational functions $F(\mathbf{x})^G$ of structural constants defined by variable matrix $\mathbf{x}=(\mathbf{x}^i_{j,k})_{i,j,k=1,2,...,m}$ is generated by the system of entries of $P(\mathbf{x})\mathbf{x}(P(\mathbf{x})^{-1}\otimes P(\mathbf{x})^{-1})$ over $F$, that is the equality \[F(\mathbf{x})^G=F(P(\mathbf{x})\mathbf{x}(P(\mathbf{x})^{-1}\otimes P(\mathbf{x})^{-1}))\] holds true.

{\bf Theorem 3.3.} The transcendence degree of $F(\mathbf{x})^G$ over $F$ equals to $m^3-m^2$ and the field extension $F(\mathbf{x})^G\subset F(\mathbf{x})$ is a pure transcendental extension.

{\bf Remark 3.1.} One of the main results (Theorem 1) of [4] states that the field extension $F\subset F(\mathbf{x})^{GL(m<F)}$ is a pure transcendental extension, which we could not get by our approach. Theorem 3.3 can be considered as a complementary result to that Theorem 1.

Now let us consider two and three dimensional algebra cases.

{\bf Example 3.l.} Two dimensional ($m=2$) case.
Let \[A=\left(
\begin{array}{cccc}
  A^{1}_{1,1} & A^{1}_{1,2}& A^{1}_{2,1}& A^{1}_{2,2} \\
  A^{2}_{1,1} & A^{2}_{1,2}& A^{2}_{2,1}& A^{2}_{2,2} \\
  \end{array}
\right)\] be the matrix of structural constants with respect to a basis. In this case at $k=1$ already $2^11!=m^1$ and therefore for the rows of $P(A)$ on can take \[\mathbf{Tr_1}(A)=(A^{1}_{1,1}+A^{2}_{2,1},A^{1}_{1,2}+A^{2}_{2,2})\ \mbox{and}\ \mathbf{Tr_2}(A)=(A^{1}_{1,1}+A^{2}_{1,2},A^{1}_{2,1}+A^{2}_{2,2})\]
and $V_0$ consists of all $A$ for which \[\det P(A)=(A^{1}_{1,1}+A^{2}_{2,1})(A^{1}_{2,1}+A^{2}_{2,2})-(A^{1}_{1,2}+A^{2}_{2,2})(A^{1}_{1,1}+A^{2}_{1,2})\neq 0.\] To see the corresponding system of generators one can
evaluate \\ $P(\mathbf{x})\mathbf{x}(P(\mathbf{x})^{-1}\otimes P(\mathbf{x})^{-1})$, where $\mathbf{x}=\left(
\begin{array}{cccc}
  \mathbf{x}^1_{1,1} & \mathbf{x}^1_{1,2}& \mathbf{x}^1_{2,1} & \mathbf{x}^1_{2,2} \\
  \mathbf{x}^2_{1,1} & \mathbf{x}^2_{1,2}& \mathbf{x}^2_{2,1} & \mathbf{x}^2_{2,2} \\
  \end{array}
\right).$

On classification problem of two dimensional algebras one can see [1,2].

{\bf Example 3.2.} Three dimensional ($m=3$) case. Let \[A=\left(
\begin{array}{ccccccccc}
  A^{1}_{1,1} & A^{1}_{1,2}& A^{1}_{1,3}& A^{1}_{2,1} & A^{1}_{2,2}& A^{1}_{2,3}& A^{1}_{3,1}& A^{1}_{3,2}& A^{1}_{3,3}\\
  A^{2}_{1,1} & A^{2}_{1,2}& A^{2}_{1,3}& A^{2}_{2,1} & A^{2}_{2,2}& A^{2}_{2,3}& A^{2}_{3,1}& A^{2}_{3,2}& A^{2}_{3,3}\\
  A^{31}_{1,1} & A^{3}_{1,2}& A^{3}_{1,3}& A^{3}_{2,1} & A^{3}_{2,2}& A^{3}_{2,3}& A^{3}_{3,1}& A^{3}_{3,2}& A^{3}_{3,3}\\
  \end{array}\right)\] be the matrix of the structural constants with respect to a basis.

In this case at $k=1$ one has $2^11! < 3^1$. At $k=2$ already $2^22!+2^11!=10 > 3^2$ and the following $10$ equalities

\[\mathbf{Tr_i}(B)\otimes \mathbf{Tr_j}(B)=\mathbf{Tr_i}(A)\otimes \mathbf{Tr_j}(A)(g^{-1})^{\otimes 2}, \ \mathbf{Tr_i}(B)B=\mathbf{Tr_i}(A)A(g^{-1})^{\otimes 2},\] where $i,j=1,2,$
 \[(B^i_{j,p}B^j_{i,q})=(A^i_{j,p}A^j_{i,q})(g^{-1})^{\otimes 2},\ (B^i_{j,p}B^j_{q,i})=(A^i_{j,p}A^j_{q,i})(g^{-1})^{\otimes 2},\] \[(B^i_{p,j}B^j_{i,q})=(A^i_{p,j}A^j_{i,q})(g^{-1})^{\otimes 2},\ (B^i_{p,j}B^j_{q,i})=(A^i_{p,j}A^j_{q,i})(g^{-1})^{\otimes 2}\] hold true.

Therefore, for example, for $Q(A)$ one can take the following matrix

 \[Q(A)= \left( \begin{array}{ccccc}
A^i_{i,1}A^j_{j,1}&A^i_{i,1}A^j_{j,2}&A^i_{i,1}A^j_{j,3}&A^i_{i,2}A^j_{j,1}&A^i_{i,2}A^j_{j,2}\\ A^i_{i,1}A^j_{1,j}&A^i_{i,1}A^j_{2,j}&A^i_{i,1}A^j_{3,j}&A^i_{i,2}A^j_{1,j}&A^i_{i,2}A^j_{2,j}\\
A^i_{1,i}A^j_{j,1}&A^i_{1,i}A^j_{j,2}&A^i_{1,i}A^j_{j,3}&A^i_{2,i}A^j_{j,1}&A^i_{2,i}A^j_{j,2}\\
A^i_{1,i}A^j_{1,j}&A^i_{1,i}A^j_{2,j}&A^i_{1,i}A^j_{3,j}&A^i_{2,i}A^j_{1,j}&A^i_{2,i}A^j_{2,j}\\
A^i_{i,j}A^j_{1,1}&A^i_{i,j}A^j_{1,2}&A^i_{i,j}A^j_{1,3}&A^i_{i,j}A^j_{2,1}&A^i_{i,j}A^j_{2,2}\\
A^i_{j,i}A^j_{1,1}&A^i_{j,i}A^j_{1,2}&A^i_{j,i}A^j_{1,3}&A^i_{j,i}A^j_{2,1}&A^i_{j,i}A^j_{2,2}\\
A^i_{j,1}A^j_{i,1}&A^i_{j,1}A^j_{i,2}&A^i_{j,1}A^j_{i,3}&A^i_{j,2}A^j_{i,1}&A^i_{j,2}A^j_{i,2}\\
A^i_{j,1}A^j_{1,i}&A^i_{j,1}A^j_{2,i}&A^i_{j,1}A^j_{3,i}&A^i_{j,2}A^j_{1,i}&A^i_{j,2}A^j_{2,i}\\
A^i_{1,j}A^j_{i,1}&A^i_{1,j}A^j_{i,2}&A^i_{1,j}A^j_{i,3}&A^i_{2,j}A^j_{i,1}&A^i_{2,j}A^j_{i,2}\\
  \end{array}\right.\]

  \[ \left. \begin{array}{cccc}
A^i_{i,2}A^j_{j,3}&A^i_{i,3}A^j_{j,1}&A^i_{i,3}A^j_{j,2}&A^i_{i,3}A^j_{j,3}\\
A^i_{i,2}A^j_{3,j}&A^i_{i,3}A^j_{1,j}&A^i_{i,3}A^j_{2,j}&A^i_{i,3}A^j_{3,j}\\
A^i_{2,i}A^j_{j,3}&A^i_{3,i}A^j_{j,1}&A^i_{3,i}A^j_{j,2}&A^i_{3,i}A^j_{j,1}\\
A^i_{2,i}A^j_{3,j}&A^i_{3,i}A^j_{1,j}&A^i_{3,i}A^j_{2,j}&A^i_{3,i}A^j_{3,j}\\
A^i_{i,j}A^j_{2,3}&A^i_{i,j}A^j_{3,1}&A^i_{i,j}A^j_{3,2}&A^i_{i,j}A^j_{3,3}\\
A^i_{j,i}A^j_{2,3}&A^i_{j,i}A^j_{3,1}&A^i_{j,i}A^j_{3,2}&A^i_{j,i}A^j_{3,3}\\
A^i_{j,2}A^j_{i,3}&A^i_{j,3}A^j_{i,1}&A^i_{j,3}A^j_{i,2}&A^i_{j,3}A^j_{i,3}\\
A^i_{j,2}A^j_{3,i}&A^i_{j,3}A^j_{1,i}&A^i_{j,3}A^j_{2,i}&A^i_{j,3}A^j_{3,i}\\
A^i_{2,j}A^j_{i,3}&A^i_{3,j}A^j_{i,1}&A^i_{3,j}A^j_{i,2}&A^i_{3,j}A^j_{i,3}\\
  \end{array}\right).\]

For $P(A)^{-1}$ one can take any $3\times 3$ size nonsingular sub-matrix of \\ $(A\otimes \mathbf{Tr_1}(A))Q(A)^{-1},$ where $(A\otimes \mathbf{Tr_1}(A))=$
\[\left( \begin{array}{ccccc}
A^1_{1,1}A^i_{i,1}&A^1_{1,1}A^i_{i,2}&A^1_{1,1}A^i_{i,3}&A^1_{1,2}A^i_{i,1}&A^1_{1,2}A^i_{i,2}\\ A^2_{1,1}A^i_{i,1}&A^2_{1,1}A^i_{i,2}&A^2_{1,1}A^i_{i,3}&A^2_{1,2}A^i_{i,1}&A^2_{1,2}A^i_{i,2}\\
A^3_{1,1}A^i_{i,1}&A^3_{1,1}A^i_{i,2}&A^3_{1,1}A^i_{i,3}&A^3_{1,2}A^i_{i,1}&A^3_{1,2}A^i_{i,2}\\
  \end{array}\right.\]
\[\left. \begin{array}{ccccc}
A^1_{1,2}A^i_{i,3}&A^1_{1,3}A^i_{i,1}&A^1_{1,3}A^i_{i,2}&A^1_{1,3}A^i_{i,3}\\
 A^2_{1,2}A^i_{i,3}&A^2_{1,3}A^i_{i,1}&A^2_{1,3}A^i_{i,2}&A^2_{1,3}A^i_{i,3}\\
A^3_{1,3}A^i_{i,3}&A^3_{1,3}A^i_{i,1}&A^3_{1,3}A^i_{i,2}&A^3_{1,3}A^i_{i,3}\\
  \end{array}\right).\]

\subsection{Commutative and anti-commutative algebra cases }

For the classification purpose instead of all $m$ dimensional algebras one can consider only such commutative or anti-commutative algebras. The commutativity (anti-commutativity) of the binary operation in terms of the corresponding matrix $A$ means $A^i_{j,k}=A^i_{k,j}$ (respectively, $A^i_{j,k}=-A^i_{k,j}$) for all $i,j,k=1,2,...,m$. So in commutative (anti-commutative) algebra case for the $V$ we consider $V=$ \[\{A\in Mat(m\times m^2;F): \ A^i_{j,k}=A^i_{k,j} (\mbox{resp.}\ A^i_{j,k}=-A^i_{k,j})\ \mbox{for all}\ i,j,k=1,2,...,m.\}\]

Note that in commutative (anti-commutative) case the dimension of $V$ is $\frac{m(m+1)}{2}$ (respectively, $\frac{m(m-1)}{2}$).

To have Theorems 2.1-2.3 for these cases one can construct a map $P: V_0\rightarrow GL(m,F)$ with property (1) in a similar way as in the general algebra case.
Consider once again equality $(3)$ and
all its possible contractions with respect to $k$ upper and $k$ lower indices.

In commutative (anti-commutative) case at $k=1$ one gets the following $1!=1$  row equality: $\mathbf{Tr_1}(B)=\mathbf{Tr_1}(A)g^{-1}=\mathbf{Tr_2}(A)g^{-1}$ as far as $A^i_{j,k}=A^i_{k,j}$ (respectively, $\mathbf{Tr_1}(B)=\mathbf{Tr_1}(A)g^{-1}=-\mathbf{Tr_2}(A)g^{-1}$ as far as $A^i_{j,k}=-A^i_{k,j}$) for all $i,j,k=1,2,...,m$.

In $k=2$ case one gets the following $2!+1!=3$ different row equalities:
\[\mathbf{Tr_1}(B)\otimes \mathbf{Tr_1}(B)=(\mathbf{Tr_1}(A)\otimes \mathbf{Tr_1}(A))(g^{-1})^{\otimes 2}
 ,\]
 \[ \mathbf{Tr_1}(B)B=\mathbf{Tr_1}(A)A(g^{-1})^{\otimes 2}, \ (B^i_{j,p}B^j_{i,q})=(A^i_{j,p}A^j_{i,q})(g^{-1})^{\otimes 2}.\]

In any $k$ case only the number of contractions of $A^{\otimes k}$ when all $k$ different  upper indices are contracted with lower indices of different $A$ is $k!$. Once again in general it is nearly clear that the corresponding resulting system of $k!$ rows depending on variable matrix $A:=\mathbf{x}=(\mathbf{x}^i_{j,k})_{i,j,k=1,2,...,m}$, where $\mathbf{x}^i_{j,k}=\mathbf{x}^i_{k,j}$ (respectively, $ \mathbf{x}^i_{j,k}=-\mathbf{x}^i_{k,j})$ for all $i,j,k=1,2,...,m$, is linear independent over $F$. But for big enough $k$ the inequality $k!\geq m^k$ holds true as well. Therefore in general for big enough $k$ it is possible to choose $m^k$ contractions (rows) among the all contractions of $\mathbf{x}^{\otimes k}$ for which the matrix $Q(\mathbf{x})$ consisting of these $m^k$ rows is nonsingular. For the matrix $Q(\mathbf{x})$ one has equality $Q(\mathbf{y})=Q(\mathbf{x})(g^{-1})^{\otimes k}$ whenever $g\in G$, $\mathbf{y}=g\mathbf{x}(g^{-1})^{\otimes 2}$.

Therefore if $P(A)^{-1}$ stands for arbitrary $m\times m$-size nonsingular sub-matrix of  $(A\otimes (\mathbf{Tr_1}(A))^{\otimes k-2})Q(A)^{-1}$ then one has the  equality $P(B)^{-1}=gP(A)^{-1}$, where $g\in G$, $B=gA(g^{-1})^{\otimes 2}$. It implies that whenever $A\in V_0=\{A\in V: \det(P(A))\det(Q(A))\neq 0\}$ the equality
$P(B)=P(A)g^{-1}$ holds true for any $g\in G$, where $B=gA(g^{-1})^{\otimes 2}$. Note that \[V_0=\{A\in V: \det(P(A))\det(Q(A))\neq 0\}\] is a $G$-invariant, open and dense subset of $V$.

Therefore we have the following results.

{\bf Theorem 3.1'.} Two commutative (anti-commutative) algebras with the matrices of structural constants $A,B\in V_0$  are the same algebras if and only if
 \[P(A)A(P(A)^{-1}\otimes P(A)^{-1})=P(B)B(P(B)^{-1}\otimes P(B)^{-1}).\]

{\bf Theorem 3.2'.} The field of $G$-invariant rational functions $F(\mathbf{x})^G$ of the structural constants presented by the matrix $\mathbf{x}=((\mathbf{x}^i_{j,k})_{i,j,k=1,2,...,m}$ of the variable commutative
(respectively, anti-commutative)
algebras, where $\mathbf{x}^i_{j,k}=\mathbf{x}^i_{k,j}$ (respectively, $\mathbf{x}^i_{j,k}=-\mathbf{x}^i_{k,j}$) for all $i,j,k=1,2,...,m$, is generated by the system of entries of $P(\mathbf{x})\mathbf{x}(P(\mathbf{x})^{-1}\otimes P(\mathbf{x})^{-1})$ over $F$, that is the equality \[F(\mathbf{x})^G=F(P(\mathbf{x})\mathbf{x}(P(\mathbf{x})^{-1}\otimes P(\mathbf{x})^{-1}))\] holds true.

{\bf Theorem 3.3'.} In commutative (anti-commutative)
 algebra case the transcendence degree of $F(\mathbf{x})^G$ over $F$ equals to $\frac{m^2(m-1)}{2}$ (respectively, $\frac{m^2(m-3)}{2}$, $m\geq 3$) and the field extension $F(\mathbf{x})^G\subset F(\mathbf{x})$ is a pure transcendental extension.

Now let us consider two dimensional commutative algebra case.

{\bf Example 3.1'.}
Let \[A=\left(
\begin{array}{cccc}
  A^{1}_{1,1} & A^{1}_{1,2}& A^{1}_{2,1}& A^{1}_{2,2} \\
  A^{21}_{1,1} & A^{2}_{1,2}& A^{2}_{2,1}& A^{2}_{2,2} \\
  \end{array}
\right)\], where $A^{i}_{1,2}= A^{i}_{2,1}$ at $i=1,2$, be the matrix of structural constants of a commutative algebra with respect to a basis. Consider $B^{\otimes 3}=g^{\otimes 3}A(g^{-1})^{\otimes 6}$ and its all contractions on $3$ upper and $3$ lower indices. Among them in particular one gets the following $6$ equalities:
\[(B^i_{\sigma(i),p}B^j_{\sigma(j),q}B^k_{\sigma(k),r})=(A^i_{\sigma(i),p}A^j_{\sigma(j),q}A^k_{\sigma(k),r})(g^{-1})^{\otimes 3}\], where $\sigma \in S_3$- the symmetric group of permutations of symbols $i,j,k$, and \[(B^i_{p,q}B^j_{r,k}B^k_{i,j})=(A^i_{p,q}A^j_{r,j}A^k_{i,j})(g^{-1})^{\otimes 3},\ (B^i_{p,q}B^j_{r,i}B^k_{k,j})=(A^i_{p,q}A^j_{r,i}A^k_{k,j})(g^{-1})^{\otimes 3}.\] So for $Q(A)$ one can take the matrix consisting of the following $8$ rows  \[(A^i_{\sigma(i),p}A^j_{\sigma(j),q}A^k_{\sigma(k),r})_{\sigma\in S_3},(A^i_{p,q}A^j_{r,k}A^k_{i,j}), (A^i_{p,q}A^j_{r,i}A^k_{k,j})\] and for the $P(A)^{-1}$ any nonsingular $2\times 2$ -size sub-matrix of
$(A\otimes Tr_1(A))Q(A)^{-1}$ provided that $\det(Q(A))\neq 0$.

On classification of three dimensional anti-commutative algebras one can see [3].

{\bf Remark 3.2.} One can consider classification of finite dimensional (general, commutative, anti-commutative) algebras with respect to subgroups of the general linear group as well. By the use of Section 3 results such classifications can be obtained for all classical subgroups of the general linear group.

{\bf Acknowledgement.} After the publication of the first variant of this paper the author was kindly informed by the author of [4] about his work on the same topic [4], which I have failed to notice before. It helped me to understand better the problems related to the classification of finite dimensional algebras.

\begin{center} References \end{center}

[1] Michel Goze and Elisabeth Remm, 2-dimensional algebras, \emph{African Journal of Mathematical Physics}, v. 10(2011),81-91.

[2] R. Dur\'{a}n D\'{\i}az et al, Classifying quadratic maps from plane to plane, \emph{Linear Algebra and its Applications} 364 (2003),pp. 1-12.

[3] L. Hern\'{a}andez Encinas, et al, Non-degenerate bilinear alternating maps $f:V\times V\rightarrow V, dim(V)=3$, over an algebraically closed field.
\emph{Linear Algebra and its Applications} 387 (2004),pp. 689-82.

[4] V. Popov, Generic Algebras: Rational Parametrization and Normal Forms, \emph{arXiv: 1411.6570v2[math.AG]}.

\end{document}